\documentclass[11pt]{article}
\usepackage[margin=1in]{geometry}
\usepackage[utf8]{inputenc}
\usepackage[T1]{fontenc}
\usepackage[english]{babel}
\usepackage{lmodern}
\usepackage{graphicx}
\usepackage{wrapfig}
\usepackage{amsmath}
\usepackage{subcaption}
\usepackage{caption}
\usepackage{paralist}
\usepackage{enumitem}
\usepackage{threeparttable}

\usepackage[linesnumbered,ruled,vlined]{algorithm2e}

\graphicspath{{./figures/}}

\usepackage[usenames]{color}
\definecolor{hellblau}{rgb}{0.2,0.4,1}
\definecolor{dunkelblau}{rgb}{0,0,0.8}
\definecolor{dunkelgruen}{rgb}{0,0.5,0}
\usepackage[
pdftex,
colorlinks,
linkcolor=dunkelblau,
urlcolor=dunkelblau,
citecolor=dunkelgruen,
bookmarks=true,
linktocpage=true,
pdfsubject={}
]{hyperref}
\usepackage{pdfpages}
\usepackage{amsmath}
\usepackage{amsthm}
\allowdisplaybreaks
\usepackage{amsfonts}

\theoremstyle{plain}

\theoremstyle{remark}

\theoremstyle{definition}

\newcommand{\Car}[1]{{\color{blue} Carol: #1}}

\begin{document}
\title
{Counting Small Cycle Double Covers}
\author{
{\sc Jorik JOOKEN\footnote{Department of Computer Science, KU Leuven Campus Kulak-Kortrijk, 8500 Kortrijk, Belgium}\;},
{\sc Ben SEAMONE\footnote{Mathematics Department, Dawson College, Montreal, QC, Canada}\;\footnote{D\'epartement d'informatique et de recherche op\'erationnelle, Universit\'e de Montr\'eal, Montreal,
QC, Canada}\;},\\[1mm]
and {\sc Carol T. ZAMFIRESCU\footnote{Department of Mathematics, Computer Science and Statistics, Ghent University, 9000 Ghent, Belgium}\;\footnote{Department of Mathematics, Babe\c{s}-Bolyai University, Cluj-Napoca, Roumania}}\;\,\footnote{E-mail addresses: jorik.jooken@kuleuven.be; bseamone@dawsoncollege.qc.ca; czamfirescu@gmail.com}}
\date{}

\maketitle
\begin{center}
\begin{minipage}{125mm}
{\bf Abstract.} A theorem due to Seyffarth states that every planar $4$-connected $n$-vertex graph has a cycle double cover (CDC) containing at most $n-1$ cycles (a ``small'' CDC). We extend this theorem by proving that, in fact, such a graph must contain linearly many small CDCs (in terms of $n$), and provide stronger results in the case of planar $4$-connected triangulations. We complement this result with constructions of planar $4$-connected graphs which contain at most polynomially many small CDCs. Thereafter we treat cubic graphs, strengthening a lemma of Hu\v{s}ek and \v{S}\'amal on the enumeration of CDCs, and, motivated by a conjecture of Bondy, give an alternative proof of the result that every planar 2-connected cubic graph on $n > 4$ vertices has a CDC of size at most~$n/2$. Our proof is much shorter and obtained by combining a decomposition based argument, which might be of independent interest, with further combinatorial insights. Some of our results are accompanied by a version thereof for CDCs containing no cycle twice.

\smallskip


\smallskip


\end{minipage}
\end{center}

\vspace{5mm}

\section{Introduction}

In a given graph, a \emph{cycle double cover} is a collection of cycles such that each of the graph's edges is contained in exactly two cycles (in this paper, a \textit{cycle} is a connected 2-regular graph). These objects are the protagonists of one of the most famous open conjectures in graph theory: \emph{Every bridgeless graph admits a cycle double cover}. Driven by different considerations---polyhedral decompositions, algorithmic graph theory, as well as matroid theory and structural graph theory---and formulating the problem in different but equivalent ways, this conjecture was published independently by Szekeres~\cite{Sz73}, Itai and Rodeh~\cite{IR78}, and Seymour~\cite{Se79} in the seventies, but has been circulating at least since the fourties. It is known as the Cycle Double Cover Conjecture (CDCC). As is frequently done, we will abbreviate ``cycle double cover'' by ``CDC''. For more on this pivotal conjecture, we refer to C.-Q.~Zhang's excellent treatise~\cite{Zh12}.

Motivated by the work of Bondy, Seyffarth, Hu\v{s}ek, and \v{S}\'amal, we shall here concern ourselves with a special class of CDCs, namely ones of small size, and their enumeration. We will also be interested in CDCs that we call \textit{true}, i.e.\ wherein no cycle appears twice. We motivate this by noting that some structurally appealing CDCs are, in fact, true: in plane 2-connected graphs which are not cycles, the CDC obtained by considering all face boundaries is true; in cubic graphs every CDC is true; etc. 

We will call a CDC of size exactly $k$ (size at most $k$; size at least $k$) a \emph{$k$-CDC}\footnote{We point out that some authors call a $k$-CDC a CDC in which every cycle is contained in one of \textit{at most} $k$ even subgraphs, and other authors call a CDC ${\cal C}$ a $k$-CDC if the cycles of ${\cal C}$ can be coloured by $k$ colours such that no pair of cycles with a common edge has the same colour. We do not use these definitions here.} (\emph{$k^-$-CDC}; \emph{$k^+$-CDC}). Unless explicitly stated otherwise, $n$ denotes the order (i.e.\ the number of vertices) of the graph that is being discussed. In particular, when we speak of a $k$-CDC and $k$ is a function in $n$, this $n$ will always denote the order of the graph in question. In our notation, an $(n-1)^-$-CDC is what Bondy called a \emph{small} CDC in~\cite{Bo90}. A graph is \emph{even} if all of its vertices have even degree. We will denote the number of $k$-CDCs of a given graph $G$ by ${\frak c}(G; k)$ and simply write ${\frak c}(k)$ when it is clear from the context to which graph we are referring. For true CDCs, we use the same notation but replace ${\frak c}$ by $\dot{{\frak c}}$. Given a graph $G$, we write $c(G)$ for the size of a CDC of $G$ of minimum size, putting $c(G) := \infty$ if no CDC exists. 


In this paper, graphs typically contain neither loops nor multiple edges, but if this is the case, it will be explicitly mentioned. All embeddings refer to embeddings on orientable surfaces. An embedded graph $G$ has vertex set $V(G)$, edge set $E(G)$, face set $F(G)$, and (orientable) genus $\gamma(G)$. The dual of $G$ will be denoted by $G^*$. A \textit{triangulation} $G$ shall be a planar graph embedded in the plane such that every face boundary is a triangle if $G$ is simple and every face boundary is a triangle or a digon if $G$ is a multigraph in which multiple edges but no loops may occur. For a separating cycle $C$ in a plane graph $G$, we denote by ${\rm Int}(C)$ (${\rm Ext}(C)$) the plane graph obtained by taking $C$ together with its interior (its exterior). For a non-negative integer $k$, we put $[k] := \{ 0, 1, \ldots, k \}$.

\section{Planar graphs}

Any CDC of a graph with maximum degree $\Delta$ contains at least $\Delta$ cycles, so if $\mathcal{G}_n$ is the family of all bridgeless graphs on $n$ vertices, we have $\max_{G \in \mathcal{G}_n} c(G) \geq n-1$. In the opposite direction, Bondy~\cite{Bo90} made in 1990 his ``Small Cycle Double Cover Conjecture'', stating that every bridgeless graph on $n$ vertices has a small CDC, i.e.\ a CDC of size at most $n - 1$. As we explain below and Seyffarth explicitly mentions in~\cite{Se92}, it is not difficult to see that the Four Colour Theorem implies that the edges of every bridgeless planar graph can be covered exactly twice by three even subgraphs (Tarsi~\cite{Ta86} notes that the same is true for Hamiltonian graphs; he gives a not-so-easy matroid theory based proof and mentions that there is a ``very easy'' graph-theoretical proof, but does not give it; we present such a proof in the final section). It then follows from a partition result of Tao (see Lemma~1 below) that every bridgeless planar graph has a CDC of size at most $3\lfloor(n-1)/2\rfloor$. This seems to be the best general upper bound for the planar case. We now restrict ourselves to planar 4-connected graphs and extend a result of Seyffarth which shows that for these graphs one can do much better.

\subsection{Planar \boldmath$4$-connected graphs}

Motivated by Bondy's aforementioned Small Cycle Double Cover Conjecture, Seyffarth gave an elegant proof of the existence of an $(n - 1)^-$-CDC for planar 4-connected graphs~\cite{Se93}. In her proof, she assumes that the Hamiltonian cycle under consideration contains two cofacial edges having a common end-vertex $v$ of degree 4 or 5. We now show that in fact the edges incident to $v$ need not be cofacial, something of no use to Seyffarth but of use to us in order to prove our first result, an enumerative version of Seyffarth's theorem. We need the following 
theorem of Tao~\cite{Ta84}; an alternative proof thereof was given by Seyffarth~\cite{Se92}.



\bigskip

\noindent \textbf{Lemma 1} (Tao, 1984 \cite{Ta84}). \emph{If $G$ is a simple planar even graph on $n$ vertices, then $E(G)$ can be partitioned into at most $\lfloor (n - 1)/2 \rfloor$ cycles.}

\bigskip

We can now state and prove our first main result relating the number of CDCs to the number of Hamiltonian cycles. 

\bigskip

\noindent \textbf{Theorem 1.} \emph{Every simple planar $4$-connected graph on $n$ vertices and with $h$ Hamiltonian cycles satisfies
$$\sum_{k \le n - 1} {\frak c}(k) \ge \frac{h}{22}.$$}

\begin{proof}
Let $G$ be a 4-connected planar graph of order $n$. It is well-known that the minimum degree of $G$ must be 4 or 5. We fix $v \in V(G)$ to be a vertex of degree 4 or 5. We now require the following claim.

\smallskip

\noindent \textsc{Claim.} \emph{Let $G$ be a planar $4$-connected graph on $n$ vertices and $v$ a vertex of degree $4$ or $5$ in $G$. Consider the following statements. For any Hamiltonian cycle ${\frak h}$ of $G$, there exists an $(n-1)^-$-CDC $\mathcal{S}$ such that (a) ${\frak h} \in \mathcal{S}$; (b) there are cycles $C^1, C^2 \in \mathcal{S}$ such that $E(C^1) \cap E(C^2)$ consists of a single edge $e$ that is incident to $v$ and the symmetric difference of $C^1$ and $C^2$ is ${\frak h}$. Then the following statements hold. (i) If $G$ has even order, then (a) holds; (ii) if $G$ has odd order, (a) or (b) holds.}

\smallskip

\noindent \emph{Proof of the Claim.} The proof of statement (i) is given by Seyffarth~\cite{Se93} (see the proof of her Theorem~4), so henceforth we assume $G$ to have odd order $n = 2k + 1$. If the two edges $e', e''$ by which ${\frak h}$ visits $v$ are cofacial, then once more we are in a situation described by Seyffarth~\cite{Se93} and we refer to her work for the details. So we may assume that $e'$ and $e''$ are not cofacial. 

Let $P$ be a plane graph which may contain multiple edges but no loops, and which contains a Hamiltonian cycle ${\frak h}$. Then the \emph{$4$-face colouring of $P$ induced by ${\frak h}$} is produced as follows. We 4-colour the faces of $P$ by properly 2-colouring the faces in the interior of ${\frak h}$ with colours 1 and 2, and properly 2-colouring the faces in the exterior of ${\frak h}$ with colours 3 and 4. This gives us a proper colouring of the faces of $P$.

Seyffarth mentions the following observation in her paper~\cite{Se93}; we recall it here, for it will be useful to us, as well. Consider the graph $G$ from the first paragraph and its Hamiltonian cycle ${\frak h}$, and the $4$-face colouring of $G$ induced by ${\frak h}$. For an integer $j$ with $2 \le j \le 4$, let $G_{1j}$ be the subgraph of $G$ consisting of the edges incident with a face of colour 1 or a face of colour $j$, but not both. Then $G_{1j}$ is an even subgraph of $G$, and every edge of $G$ lies in exactly two of $G_{12}$, $G_{13}$, $G_{14}$. We call this property ($\dagger$).

We have $G_{12} = {\frak h}$. We remark that by Lemma~1, $E(G_{13})$ and $E(G_{14})$ can each be partitioned into at most $\lfloor (n - 1)/2 \rfloor = \lfloor 2k/2 \rfloor = k$ cycles. This gives an $n^-$-CDC of $G$. So the challenge is to describe a CDC requiring one cycle less.

We first deal with the case when $v$ has degree 4. Let $v$ have neighbours $v_0, v_1, v_2, v_3$ in cyclic order. These are labelled such that ${\frak h}$ uses $v_0vv_2$. We consider the $4$-face colouring of $G$ induced by ${\frak h}$, setting the colour of the face whose boundary contains $v_2vv_3$ ($v_3vv_0$; $v_1vv_2$; $v_0vv_1$) to 1 (2; 3; 4). Remove $v$ from $G_{13}$ and add $v_1v_3$ (note that this edge was not present in $G$ since $G$ is 4-connected). This gives a planar even simple graph $G'_{13}$. By Lemma~1
$$cd(G'_{13}) \le \lfloor((n - 1) - 1)/2\rfloor = \lfloor n/2 \rfloor - 1 = k - 1,$$
where we denote by $cd(G)$ the minimum number of cycles required in a cycle decomposition of $G$. 

Subdividing once the edge $v_1v_3$ in $G'_{13}$, we obtain the graph $G''_{13}$ which satisfies $cd(G''_{13}) \le k - 1$. By ($\dagger$), $G_{12} = {\frak h}$, $G_{13}$ and $G_{14}$ cover every edge of $G$ exactly twice. Using this fact, it is not difficult to check that ${\frak h}$, $G''_{13}$ and $G_{14}$ also do so, and that this yields a CDC of $G$ containing at most $1 + k - 1 + k = 2k = n - 1$ cycles.

It remains to deal with the case when $v$ has degree 5. Let $v$ have neighbours $v_0, \ldots, v_4$ in cyclic order. Let $G'$ denote the graph obtained from $G$ by duplicating the edge $vv_3$; clearly, ${\frak h}$ is a Hamiltonian cycle in $G'$. We now consider the 4-face colouring of $G'$ induced by ${\frak h}$. Without loss of generality, the faces of $G'$ incident with $v$ are coloured as shown in
Fig.~1. Note that we can choose both the colouring of one interior face---either 1 or 2---as well as the colouring of one exterior face---either 3 or 4.

\begin{center}
\includegraphics[height=40mm]{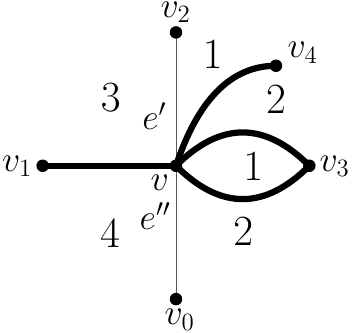}\\
Fig.~1: The vertex $v$ and its incident edges and faces. Bold edges show the graph $G'_{13}$.
\end{center}

We consider the even subgraphs $G'_{13}$ and $G'_{14}$ of $G'$. The vertex $v$ has degree 4 in $G'_{13}$, and the 2-cycle $vv_3v$ is in both $G'_{13}$ and $G'_{14}$. Let $G''_{13}$ denote the graph obtained from $G'_{13}$ by deleting the edges of the 2-cycle $vv_3v$. The degree of $v$ in $G''_{13}$ is 2. Finally, let $G'''_{13}$ be the graph obtained from $G''_{13}$ by suppressing $v$ and merging its two incident edges into one. This yields a simple graph as $v_1$ and $v_4$ are non-adjacent in $G''_{13}$ as $G$ is 4-connected (observe that $vv_1$ and $vv_4$ are not consecutive in the ordering of the edges around $v$---otherwise there could be an edge $v_1v_4$---as illustrated in Fig.~1). So by Lemma~1
$$cd(G'''_{13}) \le \lfloor((n - 1) - 1)/2\rfloor = \lfloor n/2 \rfloor - 1 = k - 1.$$

In $G'_{14}$ we subdivide one of the edges in the 2-cycle $vv_3v$ and then apply Lemma~1 to obtain
$$cd(G'_{14}) \le \lfloor ((n + 1) - 1)/2 \rfloor = \lfloor n/2 \rfloor = k.$$

Let $S_1$ be a partition of $E(G'''_{13})$ into $cd(G'''_{13})$ cycles. Then $S_1$ is a partition of the edges of $G_{13} - vv_3$ into cycles (note that we are subdividing $v_1v_4$ in order to obtain $G_{13} - vv_3$ and not $G''_{13}$). Let $S_2$ be a partition of $E(G'_{14})$ into $cd(G'_{14})$ cycles. 

If $vv_3v \notin S_2$, then $S_2$ yields a cycle cover $S'_2$ of $G_{14}$ where $vv_3$ is covered twice and all other edges are covered once. Then $S_1 \cup S'_2 \cup \{ {\frak h} \}$ is a CDC of $G$ of cardinality 
$$cd(G'''_{13}) + cd(G'_{14}) + 1 \le k - 1 + k + 1 = 2k = n - 1.$$

Now assume $vv_3v \in S_2$. As both $v$ and $v_3$ must lie in ${\frak h}$ but $vv_3$ does not, we can see ${\frak h} \cup vv_3v$ as two cycles whose intersection is $K_2$ and which together span $G$. More precisely, there exist two cycles $C^1, C^2$ such that $E(C^1) \cap E(C^2)$ consists of a single edge $e$ that is incident to $v$ and $(C^1 \cup C^2) - e = {\frak h}$. Then $S_1 \cup (S_2 \setminus \{ vv_3v \}) \cup \{ C^1, C^2 \}$ is a CDC of $G$ of size 
$$cd(G'''_{13}) + cd(G'_{14}) - 1 + 2 \le k - 1 + k - 1 + 2 = 2k = n - 1.$$

This completes the proof of the Claim.

\smallskip

For any CDC $\mathcal{C}$ of $G$, given two cycles from $\mathcal{C}$ whose union spans $G$ and whose intersection consists of precisely one edge (together with the end-vertices of this edge), at least one of these two cycles has length larger than $n / 2$. By the Claim, it suffices to show that for any fixed $t$-CDC $\mathcal{S}$ of $G$ with $t \le n - 1$, there are at most 11 cycles of length larger than $n/2$ in $\mathcal{S}$. By Euler's formula, the sum of lengths of all cycles in $\mathcal{S}$ is $2|E(G)| \le 6n - 12$. Therefore, there are at most $(6n - 12) / (n / 2) < 12$ cycles of length larger than $n / 2$ in $\mathcal{S}$. 
\end{proof} 

Our proof of Theorem~1 cannot be used for \textit{true} CDCs, at least not without significant changes. We do not know whether a true version of Theorem~1 holds, or even whether a true version of Seyffarth's theorem holds. It is clear that every planar 4-connected graph contains some true CDC, but how small (in relation to the graph's order) this can be remains open. Straightforward computer-assisted arguments show that for every $n \le 10$, every planar 4-connected graph on $n$ vertices has a true $(n-2)^{-}$-CDC. Note that due to double wheels---i.e.\ the join of a cycle and $2K_1$---which have maximum degree $n - 2$ and thus trivially no CDC with fewer than $n - 2$ cycles, no better general bound than $n - 2$ is possible. It is easy to see that double wheels do indeed admit an $(n-2)$-CDC. We note in passing that there exist (small) triangulations of maximum degree $\Delta$ in which every CDC has size greater than $\Delta$.

Recently, substantial progress has been made regarding the Hamiltonian cycle enumeration problem in planar graphs. We shall use three such results to infer some consequences of Theorem~1. Brinkmann and Van Cleemput showed that planar 4-connected graphs contain an at least linear number of Hamiltonian cycles~\cite{BV21} (while it is believed that the true bound should be quadratic), and Liu, Wang, and Yu proved that planar 4-connected triangulations contain a quadratic number of Hamiltonian cycles~\cite{LWY22}; in general, this cannot be improved since an $n$-vertex double wheel has exactly $2(n-4)(n-2)$ Hamiltonian cycles. Building on work of Alahmadi, Aldred, and Thomassen~\cite{AAT20}, it was proven by Lo and Qian~\cite{LQ22} that any planar 4-connected $n$-vertex triangulation with at most $\frac{n}{324}$ 4-cuts has $2^{\Omega(n)}$ Hamiltonian cycles. Combining these results with Theorem~1 we obtain the following. 

\bigskip

\noindent \textbf{Corollary 1.} \emph{Let ${\cal G}_n$ be the family of all $n$-vertex planar $4$-connected graphs and ${\cal T}_{n,4}$ (${\cal T}_{n,4+}$) the family of all $n$-vertex planar $4$-connected triangulations (planar $4$-connected triangulations with at most $\frac{n}{324}$ $4$-cuts). Then there is a constant $c > 0$ such that\\[1mm]
(i) each graph in ${\cal G}_n$ has at least $cn$ CDCs of size at most $n - 1$;\\
(ii) each triangulation in ${\cal T}_{n,4}$ has at least $cn^2$ CDCs of size at most $n - 1$; and\\
(iii) each triangulation in ${\cal T}_{n,4+}$ has at least $2^{cn}$ CDCs of size at most $n - 1$.}

\bigskip

Despite this progress, it remains unknown whether every planar 4-connected graph contains an $(n-2)^-$-CDC. In the remainder of this section we will treat three problems naturally arising from Theorem~1 and its consequences, namely:\\[1mm]
(I) How much larger than $n - 1$ must the size of a CDC be in order to guarantee the existence of infinitely many planar 4-connected graphs with a constant number of CDCs of that particular size;\\
(II) establish an upper bound\footnote{This upper bound is polynomial in $n$ for any fixed $c$---for quite some time, we believed that it would grow at least exponentially, in line with the typical enumerative behaviour of CDCs.} for the number of $(n + 2 - c)$-CDCs in antiprisms, where $c$ is a non-negative integer different from 2; and\\
(III) give an upper bound for the size of a CDC of minimum size of a triangulation.

\smallskip


The next theorem treats (I) as well as (II) for $c \ne 2$, while in the next subsection we address~(III). For the next theorem's proof, we will use the entropy inequality from Flum and Grohe \cite[p.~427]{FG06}, which we now recall.

\bigskip

\noindent \textbf{Lemma 2} (Flum and Grohe, 2006~\cite{FG06}). \emph{For all integers $k$ and $n$ satisfying $1 \leq k < n$ and $2k \leq n$ we have
$$\sum_{i=0}^{k} \binom{n}{i} \leq 2^{H\Big(\frac{k}{n}\Big)n}$$
where $H(p) := p \log_2{\Big(\frac{1}{p}\Big)}+(1-p) \log_2{\Big(\frac{1}{1-p}\Big)}$ is the entropy function.}

\bigskip

\noindent \textbf{Theorem 2.} (i) \textit{For all integers $k \geq 3, \ell \geq 2$, there exists a planar graph with $f :=2k(\ell - 1) + 2$ faces which satisfies $${\frak c}(f) = \frac{\ell(\ell+1)}{2} \quad {\rm \textit{and}} \quad \dot{{\frak c}}(f) = 1.$$  Furthermore, such graphs can be constructed to be $3$-connected when $k = 3$ and $\ell \geq 3$ and $4$-connected otherwise.} \\[1mm]
(ii) \textit{Let $c$ be a non-negative integer different from $2$.  There exists a graph $A$ with $|V(A)| > 6c+6$,  an exponential function $f$ in one variable $c$ (not depending on $|V(A)|$), and a function $g$ which is polynomial in $|V(A)|$ with degree at most $5c$ such that $${\frak c}(A; |V(A)|+2-c) \le f(c)\;g(|V(A)|,c).$$ Moreover, ${\frak c}(A; |V(A)|+d) = 0$ for all integers $d \ge 3$.}

\begin{proof}
We first prove part (i). Let $C_k$ be the cycle on $k$ vertices and $P_\ell$ the path on $\ell$ vertices, where $k \ge 3$ and $\ell \ge 2$. Consider the Cartesian product of $C_k$ and $P_\ell$, which we denote by $\Pi$. This graph contains pairwise disjoint $k$-cycles $C^0, \ldots, C^{\ell-1}$, of which $C^1, \ldots, C^{\ell-2}$ are separating, such that some vertex of $C^i$ is adjacent to some vertex of $C^{i+1}$ for all $i \in [\ell - 2]$. For every $i \in [\ell - 2]$ we write $C^i = ik, ik+1, \ldots, ik + k - 1, ik$ such that $ik+j$ is adjacent to $(i+1)k+j$ for every $j$. We obtain the graph $G$ by adding to $\Pi$ all edges contained in
$$E'_{\ell+i} := \{ \{ ik+j,(i+1)k+j+1 \} : j \in [k-2]\} \, \cup \, \{ \{ik+k-1,(i+1)k \} \}, \quad i \in [\ell-2].$$ The graph $G$ can be embedded in the plane by nestedly drawing $\ell$ cycles on $k$ vertices and drawing appropriate triangles in between them. Verifying the last sentence in part~(i)'s statement (on connectivity) is routine and therefore omitted; we note that, as $G$ is 3-connected, its embedding in the plane is unique~\cite{Wh33}. Fig.~2 shows an example of the graph $G$ obtained by choosing $k=4$ and $\ell=3$. For $\ell = 2$ we obtain what is known as an \textit{antiprism}, i.e.\ the square of a cycle of even length at least 6. These will play a central role in the proof of part~(ii).
\begin{center}
\includegraphics[height=40mm]{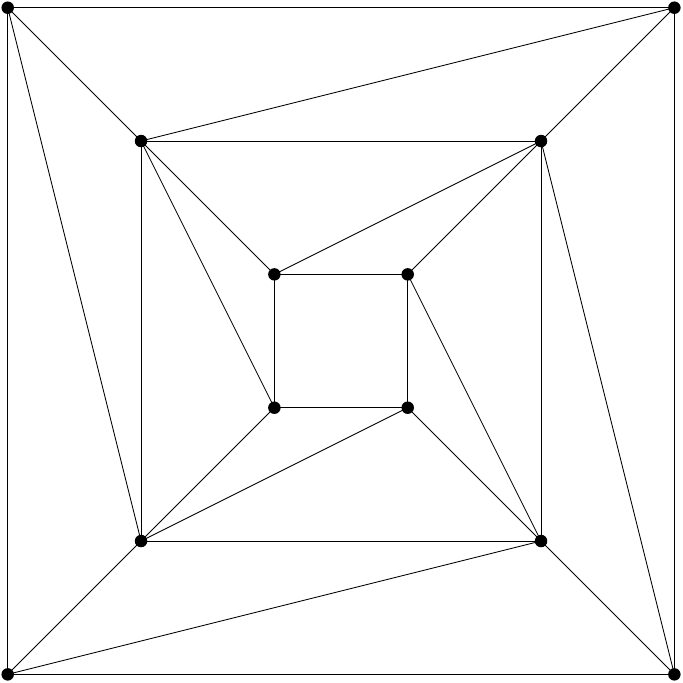}\\
Fig.~2: The graph $G$, where $k=4$ and $\ell=3$.
\end{center}
The graph $G$ has $k\ell$ vertices, $k(3\ell-2)$ edges and thus $f := 2k(\ell-1)+2$ faces by Euler's formula. Let $\mathcal{S}$ be any $f$-CDC of $G$. We write $V_i := V(C^i)$ and $E_i := E(C^i)$ for all $i \in [\ell - 1]$. Moreover, put 
$$E_{\ell + i} := E'_{\ell + i} \cup \{ \{ik+j,(i+1)k+j \} \,:\, j \in [k-1] \}.$$
We call a cycle $C \in \mathcal{S}$ {\it of type $i$}, where $i \in [\ell-1]$, iff $E(C)=E_i$ and {\it of type $\ell+i$}, where $i \in [\ell-2]$, iff $i=\min\{j \in [\ell-2] : E(C) \cap E_{\ell+j} \neq \emptyset\}$. Note that for any $i \in [\ell-2]$ and any cycle $C$ of type $\ell+i$, we have that $|E(C) \cap E_{\ell+i}|$ is an even integer that is at least 2, because $E_{\ell+i}$ is an edge cut of $G$. Let $n_i$ be the number of cycles in $\mathcal{S}$ of type $i$ with $i \in [2\ell-2]$. The following equations simultaneously hold: 
$$\sum_{i=0}^{2\ell-2}n_i=2k(\ell-1)+2,$$
$$\sum_{C \in \mathcal{S}} (|E(C) \cap E_{\ell+i}|)=4k \quad (i \in [\ell-2]), \quad {\rm and}$$
$$\sum_{C \in \mathcal{S}} (|E(C)|)=2k(3\ell-2).$$
Because of the second equation, we have $n_{\ell+i} \leq 2k$ $(i \in [\ell-2])$ and thus $\sum_{i=0}^{\ell-1}n_i \geq 2$ because of the first equation. We also have the following inequality: 
$$\sum_{C \in \mathcal{S}} (|E(C)|) \geq 3 \sum_{i=0}^{\ell-2}n_{\ell+i} + k\sum_{i=0}^{\ell-1}n_i \geq 3(\ell-1)2k + 2k = 2k(3\ell-2).$$
Equality needs to hold and this in turn means that  for every cycle $C$ of type $\ell+i$ ($i \in [\ell-2]$) we have $|E(C)|=3$ and $|E(C) \cap E_{\ell+i}|=2$ (i.e.\ $C$ is isomorphic with the triangle $K_3$). Because of the second equation, we then also have $n_{\ell+i} = 2k$ $(i \in [\ell-2])$ and $\sum_{i=0}^{\ell-1}n_i = 2$ because of the first equation. Now, for each of the $\frac{\ell(\ell+1)}{2}$ choices for choosing two cycles $C_0$ and $C_1$ such that $E(C_0)=E_i$ and $E(C_1)=E_j$ $(0 \leq i \leq j \leq \ell-1)$ we will show that there is exactly one multiset $\mathcal{S'}$ such that $C_0, C_1 \in \mathcal{S'}$ and $\mathcal{S'}$ is an $f$-CDC of $G$. Precisely one of these choices, namely when $E(C_0)=E_0$ and $E(C_1)=E_{\ell-1}$, leads to a situation where $\mathcal{S'}$ is a true $f$-CDC of $G$.

Let $T_{i,j}$ be the set of triangles that occur as a subgraph of $G$ such that each triangle in $T_{i,j}$ contains precisely $2-j$ vertices from $V_i$ and $1+j$ vertices from $V_{i+1}$, where $i \in [\ell-2]$ and $j \in [1]$. Initially we set $\mathcal{S'}=\{C_0,C_1\}$ and we will consecutively consider triangles in the sets $T_{0,0}, T_{0,1}, T_{1,0}, T_{1,1}, \ldots, T_{\ell-2,1}$ to add to $\mathcal{S'}$. Suppose we are currently considering for each triangle $T \in T_{i,0}$ how often it should occur in $\mathcal{S'}$. Any fixed $C_0$ and $C_1$ force this choice. More specifically, for each edge $e \in E_i$ with $i \in [\ell - 2]$, if 
$$\sum_{\substack{C \in \mathcal{S'}}} |E(C) \cap \{e\}|  = 2-m$$
for some $m \in [2]$, then the triangle $T \in T_{i,0}$ that contains $e$ as an edge must be added precisely $m$ times to $\mathcal{S'}$. This in turn forces that every triangle in $T_{i,1}$ must be added precisely $2-m$ times to $\mathcal{S'}$. Hence, after considering all triangles there are precisely $2k$ triangles in $\mathcal{S'} \cap (T_{i,0} \cup T_{i,1})$ for each $i \in [\ell-2]$ and $\mathcal{S'}$ is indeed an $f$-CDC of $G$ containing $C_0$ and $C_1$. Moreover, $\mathcal{S'}$ is a true $f$-CDC of $G$ when $C_0$ is different from $C_1$ and no triangle is contained twice in $\mathcal{S'}$, which happens precisely when $E(C_0)=E_0$ and $E(C_1)=E_{\ell-1}$. This completes the proof of part (i).

We now prove part (ii). Let $A$ be the graph as defined in (i) by choosing $\ell=2$ and $k > 3c+3$, i.e.\ a sufficiently large antiprism on $n := 2k$ vertices. For $c=0$, the statement follows directly from part~(i), so we henceforth assume $c \geq 1$. We define $V_0, V_1, E_0, E_1$ and $E_2$ as in the proof of part~(i). We call a cycle $C$ of $A$ a \textit{rare cycle} if $|E(C)| \geq k$ and at least one of the cycles $A[V_0]$ and $A[V_1]$ lie in the interior of $C$ in the embedding of $A$ that was discussed in part~(i). If $A$ does not have an $(n+2-c)$-CDC, then the statement is trivially true. Otherwise, let $\mathcal{S}$ be an $(n+2-c)$-CDC of $A$. Since $\mathcal{S}$ is a CDC, the following equations hold:
\begin{gather}
\sum_{C \in \mathcal{S}} (|E(C) \cap E_{0}|)=2k \label{eq:E0} \\
\sum_{C \in \mathcal{S}} (|E(C) \cap E_{2}|)=4k \label{eq:E2} \\
\sum_{C \in \mathcal{S}} (|E(C)|)=8k \label{eq:all}
\end{gather}
Let $m(C,\mathcal{S})$ represent how often the cycle $C$ occurs in $\mathcal{S}$ (i.e.\ 0, 1 or 2 times), let\\$s = \sum_{C \in \mathcal{S}\text{ and }|E(C)| \geq k} (m(C,\mathcal{S}))$ and let $t = \sum_{C \in \mathcal{S}\text{ and }|E(C)| < k} (|E(C)|)$. Now we have $t \leq 8k-sk$ because of Equation~(\ref{eq:all}), but also $t \geq 3(2k+2-c-s)$ because any cycle has at least 3 edges. If $s \geq 3$ would hold, then we would obtain a contradiction, because 
$$k>3c+3 \Rightarrow k > \frac{3c+3(s-2)}{s-2} \Rightarrow sk-8k > -3(2k+2-c-s) \Rightarrow 8k-sk < 3(2k+2-c-s).$$
Therefore, we have $s \leq 2$.

Note that for every cycle $C$ of $A$, we have that $|E(C) \cap E_2|$ is an even integer, because $E_2$ is an edge cut of $A$. Furthermore, if $|E(C) \cap E_2| =0$, then $|E(C)|=k$ and if $|E(C) \cap E_2| \geq 4$, then $|E(C)| \geq k+1$ (in both cases $C$ is a rare cycle). Let $s' = \sum_{C \in \mathcal{S}\text{ and }C\text{ is rare}} (m(C,\mathcal{S}))$. Since $s \leq 2$, we also have $s' \leq 2$. We will now show that $s'=2$, because $s' \leq 1$ would lead to a contradiction. 

Assume $s' \leq 1$. If $c=1$, then let $C^{0} = A[V_0]$ and $C^{1} = A[V_1]$ be the only two cycles in $A$ that do not contain any edge from $E_2$. It is impossible that $m(C^{0},\mathcal{S})+m(C^{1},\mathcal{S})=0$, because of Equation~(\ref{eq:E2}). Hence, we assume without loss of generality that $m(C^{0},\mathcal{S})=1$ and $m(C^{1},\mathcal{S})=0$. Note that $|E(C^{0}) \cap E_{0}| = k$ and that the $2k$ other cycles in $\mathcal{S}$ each contain exactly 2 edges from $E_2$ because of Equation~(\ref{eq:E2}). Let 
$$T := \{C : C\text{ is a triangle in }A\text{ for which }|E(C) \cap E_1|=1\text{ and }|E(C) \cap E_0|=0\}$$
be the set of the only cycles in $A$ that are not rare and do not contain any edge from $E_0$. If $\sum_{C \in T} m(C,\mathcal{S}) < k$, we obtain the following contradiction with Equation~(\ref{eq:E0}):
$$\sum_{C \in \mathcal{S}} (|E(C) \cap E_{0}|) \geq k+(2k-\sum_{C \in T} m(C,\mathcal{S})) >2k$$
Hence, $\sum_{C \in T} m(C,\mathcal{S}) \geq k$. We now consider two cases. 

Case 1: There is a cycle $C \in T$ for which $m(C,\mathcal{S})=2$. Let $e_0 \in E_1$ and $e_1, e_2 \in E_2$ be the three distinct edges of $C$ and $v \in V_0$ be a vertex of $C$. There must be another cycle $C' \in \mathcal{S}$ containing the two edges $e_3, e_4 \in E_0$ which are incident with $v$ ($e_1$, $e_2$, $e_3$ and $e_4$ are pairwise distinct), but $C'$ cannot contain $e_0$. Therefore at least one of the cycles $A[V_0]$ and $A[V_1]$ lie in the interior of $C'$ and $|E(C')| \geq k$. Hence $C^{0}, C' \in \mathcal{S}$ must be rare cycles and this gives a contradiction. 

Case 2: There is no cycle $C \in T$ for which $m(C,\mathcal{S})=2$. Since $\sum_{C \in T} m(C,\mathcal{S}) \geq k$ and $|T|=k$, we must have $m(C,\mathcal{S})=1$ for each cycle $C \in T$. Define $\mathcal{O} := \{C : C \in \mathcal{S}\text{ and }C \neq C^{0}\text{ and } C \notin T\}$ to be cycles in $\mathcal{S}$ that we did not yet consider. Note that every edge $e \in E(A)$ is contained in exactly one cycle in $\mathcal{O}$. If $\mathcal{O}$ contains a triangle, we obtain a contradiction by a very similar argument as the one that was used in Case~1. If $C \in \mathcal{O}$ is a cycle of length at least 4, then every edge from the set $E(A[V(C)]) \setminus E(C)$ must be contained in the same cycle $C' \in \mathcal{O}$. But then at least one of the cycles $A[V_0]$ and $A[V_1]$ lie in the interior of $C'$ and $|E(C')| \geq k$. Hence $C^{0}, C' \in \mathcal{S}$ must be rare cycles and this gives a contradiction.

If $c>2$, then $2(2k+2-c) < 4k$, so there must be some cycle $C \in \mathcal{S}$ for which $|E(C) \cap E_2| \geq 4$ because of Equation~(\ref{eq:E2}) and thus $C$ must be a rare cycle. Now there are four cases: $C$ contains three consecutive edges $e_1$, $e_2$ and $e_3$ such that 1) $e_1 \in E_0, e_2 \in E_2$ and $e_3 \in E_1$, 2) $e_1 \in E_0, e_2 \in E_2$ and $e_3 \in E_2$, 3)  $e_1 \in E_1, e_2 \in E_2$ and $e_3 \in E_0$ and 4) $e_1 \in E_1, e_2 \in E_2$ and $e_3 \in E_2$. Case 3 is symmetric with Case 1, and Case 4 is symmetric with Case 2. In what follows, we deal with Case~1 and Case 2 (see Fig.~3a and Fig.~3b, respectively). Each of these cases has several subcases, but the argument for all of the subcases is nearly the same, except for some variables $a_i$ that need to be changed. The argument is as follows: The cycle $C$ contains $e_1, e_2$ and $e_3$. Consider the second cycle $C'$ that contains $e_2$ (note that $C'$ is not rare). If $C'$ contains the edges $a_1$ and $a_2$, then there must be two distinct cycles $C_1$ and $C_2$ such that $C_1$ contains the edges $a_3$ and $a_4$, whereas $C_2$ contains the edges $a_3$ and $a_5$, but only one of them can also contain the edge $a_6$, so the other one must be a rare cycle and we obtain a contradiction. All possible values for $a_i$ that cover all the subcases are given in Table~\ref{tab:subcases}.

\begin{figure}[h]
\centering
\begin{subfigure}{0.3\linewidth}
\centering\includegraphics[width=1.1\linewidth]{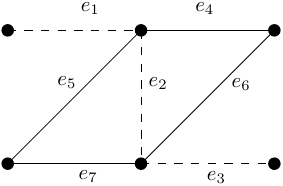}
\caption{}
\label{fig:fragment1}
\end{subfigure}%
\hspace{2cm}
\begin{subfigure}{0.3\linewidth}
\centering\includegraphics[width=1.1\linewidth]{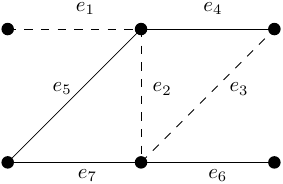}
\caption{}
\label{fig:fragment2}
\end{subfigure}
\label{fig:fragments}
\centering \\Fig.~3: Two cases for a cycle $C \in \mathcal{S}$ having $|E(C) \cap E_2| \geq 4$. The horizontal edges on the top are in $E_0$, the horizontal edges on the bottom are in $E_1$ and the others are in $E_2$.
\end{figure}

\begin{table}[h!] \centering
	\begin{threeparttable}
		\begin{tabular}{cccccccccccccccccccccc} \\
			\hline
			\noalign{\smallskip}
			Case & Figure & $a_1$ & $a_2$ & $a_3$ & $a_4$ & $a_5$ & $a_6$ \\
			\noalign{\smallskip}
			\hline
			\noalign{\smallskip}
			\multicolumn{1}{c}{Case 1a} & Fig. 3a & $e_3$ & $e_4$ & $e_7$ & $e_6$ & $e_6$ & $e_4$\\
			\multicolumn{1}{c}{Case 1b} & Fig. 3a & $e_6$ & $e_4$ & $e_7$ & $e_6$ & $e_3$ & $e_4$\\
			\multicolumn{1}{c}{Case 1c} & Fig. 3a & $e_7$ & $e_1$ & $e_4$ & $e_5$ & $e_5$ & $e_7$\\
			\multicolumn{1}{c}{Case 1d} & Fig. 3a & $e_7$ & $e_5$ & $e_4$ & $e_5$ & $e_1$ & $e_7$\\
			\multicolumn{1}{c}{Case 2a} & Fig. 3b & $e_3$ & $e_4$ & $e_7$ & $e_6$ & $e_6$ & $e_4$\\
			\multicolumn{1}{c}{Case 2b} & Fig. 3b & $e_6$ & $e_4$ & $e_7$ & $e_6$ & $e_3$ & $e_4$\\
			\multicolumn{1}{c}{Case 2c} & Fig. 3b & $e_7$ & $e_1$ & $e_4$ & $e_5$ & $e_5$ & $e_7$\\
			\multicolumn{1}{c}{Case 2d} & Fig. 3b & $e_7$ & $e_5$ & $e_4$ & $e_5$ & $e_1$ & $e_7$\\
			\hline
		\end{tabular}
	\end{threeparttable}
	\caption{A summary of the subcases.}	
	\label{tab:subcases}
\end{table}

Hence, we conclude that $s'=2$. The intuition for the rest of the proof is as follows: we will first derive an upper bound for the number of ways to choose the two rare cycles, then we will derive an upper bound for the number of ways to choose all cycles with length strictly larger than 3 and after these two choices we will see that there is only at most one way to choose cycles with length precisely 3. The product of these upper bounds is then an upper bound for the number of $(n+2-c)$-CDCs of $A$.

Let $C$ and $C'$ be the two (not necessarily distinct) rare cycles in $\mathcal{S}$. Note that $|E(C) \cap E_2|+|E(C') \cap E_2| \leq 2c$ because of Equation~\ref{eq:E2}. For every set of edges $S \subset E_2$, there are at most two distinct rare cycles $C''$ and $C'''$ such that $E(C'') \cap E_2 = S$ and $E(C''') \cap E_2 = S$. By iterating over all values for $l_1 := |E(C) \cap E_2|$ and $l_2 := |E(C') \cap E_2|$ whose sum $s$ is at most $2c$, we get an upper bound for the number of ways to choose two rare cycles $C, C' \in \mathcal{S}$:

\begin{equation*}
\begin{aligned}
\sum_{\substack{l_1=0 \\ l_1 \text{ even}}}^{2c} \sum_{\substack{l_2=l_1 \\ l_2 \text{ even}}}^{2c-l_1} 2^2 \binom{2k}{l_1} \binom{2k}{l_2} &\leq 4 \sum_{s=0}^{2c} \sum_{l_1=0}^{s} \binom{2k}{l_1} \binom{2k}{s-l_1}\\ &\overset{\text{(1)}}{=} 4 \sum_{s=0}^{2c} \binom{4k}{s} \\
&\overset{\text{(2)}}{\leq} 4 \cdot 2^{(\frac{2c}{4k}\log_2(\frac{4k}{2c}) + \frac{4k-2c}{4k}\log_2(\frac{4k}{4k-2c}))4k}\\
&= 4 \Big(\frac{2k}{c}\Big)^{2c}\Big(1+\frac{2c}{4k-2c}\Big)^{4k-2c}\\
&\overset{\text{(3)}}{\leq} 4 \Big(\frac{2k}{c}\Big)^{2c} e^{2c}\\
&= 4 \Big(\frac{e}{c}\Big)^{2c} n^{2c},
\end{aligned}
\end{equation*}

where we used (1) Vandermonde's identity, (2) the entropy inequality from Lemma 2, and (3) the fact that $(1+\frac{m}{n})^{n} \leq e^{m}$. 

For each fixed length $q$, there are at most $2k$ cycles in $A$ which are not rare and have length exactly $q$. Hence, if the lengths of the cycles in $\mathcal{S}$ are fixed, it is easy to obtain an upper bound for the number of ways to choose cycles in $\mathcal{S}$ which are not rare, but have length strictly greater than~3. Let $q_1 \geq q_2 \geq \ldots \geq q_{2k+2-c}$ be the lengths of the cycles in $\mathcal{S}$.
Since $C$ and $C'$ are rare cycles, we have $2k \leq q_1+q_2$. The sum of the lengths of the other $2k-c$ cycles in $\mathcal{S}$ is bounded from below as follows:
$$\sum_{i=3}^{2k+2-c} q_i \geq 6k-3c = 3(2k-c)$$
because every cycle has length at least 3. Define $q_{i}' := q_i-3$ (for all $i \in \{3,4,\ldots,2k+2-c\}$) and $r := \sum_{i=3}^{2k+2-c} q_{i}'$. Now we have $3c \geq r \geq 0$ because of Equation~(\ref{eq:all}). By iterating over all possible values that $r$ could have, we obtain an upper bound for the number of ways to choose cycles in $\mathcal{S}$ which are not rare, but have length strictly greater than 3:
$$\sum_{r'=0}^{3c} p(r')(2k)^{r'} \leq (3c+1)p(3c)n^{3c}.$$
Here, $p(r')$ is the partition function that represents the number of ways to write $r'$ as a sum of strictly positive integers, which is also equal to the number of ways to write $r'$ as a sum of $2k-c$ non-negative integers $q_{i}'$ ($i \in \{3,4,\ldots,2k+2-c\}$), because $2k-c > r'$.

Finally, the only cycles that are left to consider are the $2k$ triangles of $A$. Their choice is forced. Let $T$ be a triangle of $A$ and let $e \in E(T)$ be an edge such that $e$ is not contained in any other triangle. Then we have $m(T,\mathcal{S})=2-j$, where $j$ is the number of times that $e$ was already covered by a cycle in $\mathcal{S}$ which was not a triangle. Hence, there is at most one way to choose the triangles in $\mathcal{S}$ after fixing all other cycles.

We conclude that the number of $(2k+2-c)$-CDCs of $A$ is bounded from above by the product of these upper bounds:
$$4 \Big(\frac{e}{c}\Big)^{2c} n^{2c} (3c+1)p(3c)n^{3c} = 4 \Big(\frac{e}{c}\Big)^{2c} (3c+1)p(3c) n^{5c}.$$
The theorem now follows by choosing $f(c) =  4 \Big(\frac{e}{c}\Big)^{2c} (3c+1)p(3c)$ and $g(n,c)=n^{5c}$.

We conclude by proving the theorem's last statement, namely that an antiprism $A$ on $2k$ vertices has no $(2k+3)^+$-CDC. Suppose $\mathcal{S}$ is a CDC of $A$ containing at least $2k+3$ cycles. Let $E_2$ be defined as in the proof of part (i). Since $E_2$ is an edge-cut, every cycle of $A$ contains an even number of edges from $E_2$. By Equation~(\ref{eq:E2}), at least three cycles of $\mathcal{S}$ should not contain any edges of $E_2$. The sum of the lengths of these cycles is at least $3k$. But then the sum of the lengths of all cycles in $\mathcal{S}$ would exceed $8k$, since every cycle has length at least 3. This contradicts Equation~(\ref{eq:all}) and finishes the proof of part (ii).\end{proof}

For statement (ii) of Theorem~2, we were not able to give a better upper bound by treating true CDCs instead of CDCs. We now state two direct corollaries of Theorem~2.

\bigskip

\noindent \textbf{Corollary 2.} \textit{For every integer $\ell \ge 2$ there exists for infinitely many $n$ a planar $4$-connected graph on $n$ vertices with exactly $\ell(\ell+1)/2$ CDCs of size $(2\ell-2)n/\ell + 2$. In particular, for infinitely many~$n$ there is a planar $4$-connected graph on $n$ vertices with exactly three $(n+2)$-CDCs, namely the antiprism. If the triangles of an $n$-vertex embedded antiprism $A$, in their cyclic order, are denoted $T_1, \ldots, T_n$, the inner $n/2$-gon $I$, and the outer $n/2$-gon $O$, and assuming that $T_1$ and $I$ share an edge, then the three $(n+2)$-CDCs of $A$ are given by $$\{ I, T_1, \ldots, T_n, O \}, \ \{ T_i, T_i \}_{i \; {\rm odd}} \cup \{ O \}, \ \{ T_i, T_i \}_{i \; {\rm even}} \cup \{ I \},$$
the latter two being equivalent under a suitable automorphism. In particular, antiprisms have a unique true $(n+2)$-CDC, where $n + 2$ coincides with the number of faces. Moreover, all graphs described in Theorem~2~(i) have a unique true $f$-CDC, where $f$ is the number of faces.} 

\bigskip

In Corollary 2 we give an answer to (I), but it is unclear whether $n + 2$ is optimal. Computational experiments, executed (and indeed executable) only for small orders, seem to indicate that there might be an infinite family of planar 4-connected $n$-vertex graphs with a constant number of $(n+1)$-CDCs but perhaps none with a constant number of $n$-CDCs. We end this subsection by mentioning another consequence of Theorem~2.


\bigskip

\noindent \textbf{Corollary 3.} \textit{Antiprisms have $O(n^{15})$ $(n-1)$-CDCs.}

\bigskip

In Theorem~2~(ii) we did not treat the case $c = 2$, i.e.\ $n$-CDCs. For this case the proof above breaks down and we do not know of an alternative route.

\subsection{Triangulations}

Triangulations may contain no $(n-2)^-$-CDC (consider a triangulation of maximum degree $n-1$), but Bondy and Seyffarth~\cite{Bo90} gave a short proof of the fact that an $n$-vertex triangulation of any surface admits an $(n-1)^-$-CDC. It follows from a result of Seyffarth~\cite{Se89} that if $G$ is a triangulation with maximum degree $\Delta \ge 8$ and of order $n > \frac{3\Delta}{2} + 1$, then $G$ admits an $(n - 2)^-$-CDC. Equivalently, $n$-vertex triangulations with maximum degree at most $(2n/3) - 1$ admit an $(n - 2)^-$-CDC. We here follow a different route to bound the size of a smallest CDC in a given triangulation, namely via Jackson-Yu decomposition trees~\cite{JY02}, a way to capture the global structure of a triangulation which we now recall. 
 


Given a triangulation $T$, we define its (unique) \textit{decomposition tree} $D$ as follows. We successively split $T$ along separating triangles, where a copy of the separating triangle is retained in each part. The vertices of $D$ correspond to the \textit{pieces} we obtain after this splitting procedure. Each of these pieces is a triangulation, either 4-connected or isomorphic to $K_4$. At every step of the splitting procedure, we split a triangulation into two parts $T_1$ and $T_2$ along a separating triangle $\Delta$, which Jackson and Yu call the \textit{marker triangle} in $T_1$ and $T_2$. When two pieces share a marker triangle we connect the corresponding vertices in $D$. For more details we refer to~\cite{JY02}. We denote by ${\cal P}(T)$ (${\cal P}_4(T)$; ${\cal P}_{K_4}(T)$) the set of all pieces (all 4-connected pieces; all pieces isomorphic to $K_4$).

\bigskip

\noindent \textbf{Theorem 3.} \textit{For a triangulation $T$, we have $$c(T) \le \sum_{G \in {\cal P}(T)} c(G).$$ 
Under the assumption that every planar $4$-connected $n$-vertex triangulation admits a CDC of size at most $n - 2$, we have
$$c(T) \le 3\,|{\cal P}_{K_4}(T)| + \sum_{G \in {\cal P}_4(T)} (|V(G)| - 2).$$}

\begin{proof}

We begin by proving the first inequality. If the plane triangulation $T$ is 4-connected there is nothing to show, so let $T$ contain a separating triangle $S_\Delta$. Then there exist plane triangulations $T_1$ and $T_2$, each on at least four vertices, each containing a facial triangle $\Delta_1$ and $\Delta_2$, respectively, such that taking the disjoint union of $T_1$ and $T_2$ and identifying $\Delta_1$ and $\Delta_2$ yields $T$. For the moment, we shall treat $T_1$ and $T_2$ as disjoint (i.e.\ prior to identification). 

In this paragraph, assume $i \in \{ 1, 2 \}$. Let ${\cal C}_i$ be a CDC of $T_i$. Let the edges of $\Delta_i$ be $e_i, f_i, g_i$ such that, when performing the identification, $e_1$ is identified with $e_2$, $f_1$ is identified with $f_2$, and $g_1$ is identified with $g_2$. Let ${\cal C}:={\cal C}_1\cup{\cal C}_2$. We now perform the identification prescribed above, but see $T_1$ and $T_2$ as subtriangulations of $T$, with $\Delta_1 = \Delta_2 = S_\Delta$, and denote by $e$ ($f$; $g$) the edge resulting from the identification of $e_1$ and $e_2$ ($f_1$ and $f_2$; $g_1$ and $g_2$). ${\cal C}$ covers every edge of $T - S_\Delta$ twice and every edge in $S_\Delta$ four times. Put $|{\cal C}_i| =: c_i$.

If the triangle $S_\Delta$ occurs at least twice in ${\cal C}$, we remove these two triangles from ${\cal C}$ and have obtained a CDC of $T$ of size $c_1 + c_2 - 2$. Henceforth, we assume the triangle $S_\Delta$ to occur at most once in ${\cal C}$. We now have two cases.


Case 1: There exist $C^1 \in \mathcal{C}_1, C^2 \in \mathcal{C}_2$ such that $|E(C^1) \cap E(C^2)|=2$. We may assume without loss of generality that $E(C^1) \cap E(C^2)=\{e,f\}$. We claim that there exist cycles $C^3 \in \mathcal{C}_1, C^3 \neq C^1$ and $C^4 \in \mathcal{C}_2, C^4 \neq C^2$ such that $E(C^3) \cap E(C^4) = \{g\}$. Since the triangle $S_\Delta$ occurs at most once in $\mathcal{C}$, there exists a $j \in \{1,2\}$ such that $\mathcal{C}_j$ contains two distinct cycles $C^5$ and $C^6$ (both different from $C^1$) which contain $g$. We may assume without loss of generality that $j=1$. Let $C^7$ be a cycle in $\mathcal{C}_{2}$ different from $C^2$ which contains $g$, but not simultaneously $e$ and $f$. If $C^7$ contains neither $e$ nor $f$, then $E(C^5) \cap E(C^7)=\{g\}$. If $C^7$ contains $e$, but not $f$, either $C^5$ or $C^6$ does not contain $e$ and thus either $E(C^5) \cap E(C^7)=\{g\}$ or $E(C^6) \cap E(C^7)=\{g\}$. Similarly, if $C^7$ contains $f$, but not $e$, either $C^5$ or $C^6$ does not contain $f$ and thus either $E(C^5) \cap E(C^7)=\{g\}$ or $E(C^6) \cap E(C^7)=\{g\}$, which proves the claim. Now the graph $C^1 \Delta C^2$ is either a cycle (in which case we define $S_0 := \{C^1 \Delta C^2\}$) or the union of two cycles $C^8$ and $C^9$ (in which case we define $S_0 := \{C^8, C^9\}$). Similarly, the graph $C^3 \Delta C^4$ is either a cycle (in which case we define $S_1 := \{C^3 \Delta C^4\}$) or the union of two cycles $C^{10}$ and $C^{11}$ (in which case we define $S_1 := \{C^{10}, C^{11}\}$), so
$$\left( {\cal C} \setminus \{ C^1, C^2, C^3, C^4 \} \right) \cup S_0 \cup S_1$$
is a CDC of $T$ of size at most $c_1 + c_2 - 4 + 2 + 2 = c_1 + c_2$.

Case 2: There do not exist $C^1 \in \mathcal{C}_1, C^2 \in \mathcal{C}_2$ such that $|E(C^1) \cap E(C^2)|=2$. This means that $|E(C^1) \cap E(C^2)| \leq 1$ for any $C^1 \in \mathcal{C}_1, C^2 \in \mathcal{C}_2$. Let $C^3, C^4, C^5 \in \mathcal{C}_1$ be three pairwise distinct cycles that contain respectively $e, f$ and $g$. Similarly, let $C^6, C^7, C^8 \in \mathcal{C}_2$ be three pairwise distinct cycles that contain respectively $e, f$ and $g$. Now, for each $j \in \{0,1,2\}$ the graph $C^{3+j} \Delta C^{6+j}$ is either a cycle (in which case we define $S_j := \{C^{3+j} \Delta C^{6+j}\}$) or the union of two cycles $C^{9+j}$ and $C^{12+j}$ (in which case we define $S_j := \{C^{9+j}, C^{12+j}\}$). Hence,
$$\left( {\cal C} \setminus \{ C^3, C^4, C^5, C^6, C^7, C^8 \} \right) \cup S_0 \cup S_1 \cup S_2$$
is a CDC of $T$ of size at most $c_1 + c_2 - 6 + 2 + 2 + 2 = c_1+c_2$.

The theorem's second inequality follows from the first inequality, the imposed assumption that planar 4-connected triangulations on $n$ vertices admit a CDC of size at most $n - 2$, and the fact that $c(K_4) = 3$.
\end{proof}

The strength of Theorem~3 very much depends on how well we can cover each piece. We have formulated the second part of the theorem as we strongly believe that the $n - 2$ bound holds---in fact, it seems very likely that it holds even in the larger class of planar 4-connected graphs. If this is the case and the triangulation contains few pieces isomorphic to $K_4$, the bound we obtain through Theorem~3 is much better than the bound one obtains by applying directly the Bondy-Seyffarth result to the entire triangulation, i.e.\ order minus 1. Note that applying the Bondy-Seyffarth result to the pieces and then applying Theorem~3 is not advantageous. To illustrate this, consider a triangulation $T$ with the simplest non-trivial decomposition tree, $K_2$. Let the order of the pieces be $n_1$ and $n_2$. For each piece the Bondy-Seyffarth bound gives us a CDC containing at most order minus 1 cycles. By Theorem~3 the entire triangulation admits a CDC of size at most $n_1 - 1 + n_2 - 1 = n_1 + n_2 - 2 = |V(T)| + 1$, but applying Bondy-Seyffarth directly to $T$ yields a CDC of size at most $|V(T)| - 1$. 

A version of Theorem~3 for true CDCs also holds, with a proof analogous to the one given above, under the assumption that each piece is covered by a true CDC.

\section{Cubic graphs}

\subsection{Bondy's conjecture on the size of CDCs of \boldmath{2}-connected cubic graphs}

Bondy conjectured that every simple $2$-connected cubic graph that is not $K_4$ has a CDC of size at most~$n/2$, see~\cite[Conjecture 5.1]{Bo90}; we will abbreviate this conjecture by (${\frak B}$). We recall that, clearly, every CDC in a cubic graph is true. 
Lai, Yu, Zhang~\cite{LYZ94} present two short proofs---one due to them, one due to Goddyn and Richter---of the fact that every CDC of a cubic $n$-vertex graph contains at most $\frac{n}{2} + 2$ cycles. 
In their paper, Lai, Yu, and Zhang then present their main result stating that if a $2$-connected cubic graph has a CDC, then it has a CDC of size at most $n/2$. 
We recall the following important conjecture, where an embedding (of a graph) is \textit{strong} if each face of the embedding is a disk whose boundary is a cycle in the graph.

\bigskip

\noindent \textbf{Strong Embedding Conjecture.} \emph{Every $2$-connected graph has a strong embedding into some orientable surface.}

\bigskip

If we assume the Strong Embedding Conjecture to be true, then the results mentioned in the first paragraph of this section imply a positive answer to Bondy's original conjecture (${\frak B}$). We now present a proposition with a short proof---based on a similar, old idea which was recently used in~\cite{HS24}---which (i) gives an alternative approach to settling (${\frak B}$) for the planar case, detailed in Theorem~4, and (ii) yields a strengthening of~\cite[Corollary~20]{HS24}, see Proposition~2.

\bigskip

\noindent \textbf{Proposition 1.} \emph{Let $G$ be a $2$-connected cubic embedded graph on $n$ vertices and $G^*$ its dual. If the Strong Embedding Conjecture holds and $G^*$ contains induced wheels $W_e^1, \dots, W_e^k$ of even order and induced cycles $C_e^1, \dots, C_e^m$ of even order, all of which are pairwise disjoint, then $G$ has a CDC of size at most
$$\frac{n}{2} - 2\gamma(G) + 2 + 3k + 2m - \sum_{i=1}^k |V(W_e^i)| - \sum_{j=1}^m |V(C_e^j)|.$$
Furthermore, $G$ has at least $2^{k + m}$ CDCs of size at most $\frac{n}{2} - 2\gamma(G) + 2$.}

\begin{proof}
By Euler's formula, $G$ has $f := \frac{n}{2} - 2\gamma(G) + 2$ faces, all of which are bounded by cycles as we are assuming the Strong Embedding Conjecture to hold. Denote by ${\cal F}$ the collection of these boundary cycles, which forms a CDC of $G$. Its size is $f$. We require the following definitions.

Let $H$ be a cubic graph, $k \ge 3$ an integer, and $C^1, \dots, C^k$ cycles in $H$ such that the intersection of $C^i$ and $C^{i + 1}$ is isomorphic to $K_2$ and the intersection of $C^i$ and $C^j$ is empty for any $i \in \{1, \dots, k\}$ and any $j \in [k + 1] \setminus \{i - 1, i, i + 1\}$, where $C^0 := C^k$ and $C^{k + 1} := C^1$. We call $\{C^1, \dots, C^k\}$ a \emph{cycle ring} and $k$ its \emph{length}. For any $i \in \{1, \dots, k - 1\}$ we denote the edge lying on both $C^i$ and $C^{i+1}$ by $e_i$, and the edge lying on $C^k$ and $C^1$ shall be $e_k$. Removing from $\bigcup_i C^i$ all edges $e_1, \dots, e_k$, we obtain two disjoint cycles, one of which (chosen arbitrarily) we call $B$. We call $\{B, C^1, \dots, C^k\}$ a \emph{wheel ring} and $k+1$ its \emph{length}.

If $k$ is even, put $\rho[C^1, \dots, C^k] := \{D_1, D_2\}$, where $D_i$ is the symmetric difference of $B$ and $C^{i} \cup C^{i + 2} \cup \dots \cup C^{k + i - 2}$ for $i \in \{1, 2\}$. If $k$ is odd, put $\omega_C[B, C^1, \dots, C^k] := \{ O_1, O_2, O_3 \}$, where $C$ is a cycle among $C^1, \dots, C^k$, chosen arbitrarily; $O_1$ is the symmetric difference of $B$ and $C$; and $O_i$ is the symmetric difference of $B$ and $(C^{i-1} \cup C^{i+1} \cup \dots \cup C^{k + 2  - i }) \setminus C$ for $i \in \{ 2,3 \}$. We call the operation $\rho$ a \emph{cycle ring exchange} and the operation $\omega_C$ a \emph{wheel ring exchange}. Note that $D_1 \ne D_2$ and $O_1, O_2, O_3$ are pairwise distinct.

We now alter ${\cal F}$ in order to obtain on one hand CDCs of smaller size, and on the other hand many CDCs of size at most $f$. Note that all of these are CDCs of size at most $\frac{n}{2}$ if the graph is non-planar. Every induced wheel (induced cycle) in $G^*$ corresponds to a wheel ring (cycle ring) in $G$, composed of cycles residing in ${\cal F}$. We perform wheel ring exchanges on all wheel rings corresponding to $W_e^1, \dots, W_e^k$ as well as cycle ring exchanges on all cycle rings corresponding to $C_e^1, \dots, C_e^m$. For each wheel $W_e^i$, the size of our CDC decreases by $|V(W_e^i)| - 3 \ge 1$, for each cycle $C_e^j$, the size of our CDC decreases by $|V(C_e^j)| - 2 \ge 2$. Thus, the first statement is proven.

For the second statement it suffices to observe that exchanging distinct wheel ring sets and cycle ring sets yields distinct CDCs.
\end{proof}

In Proposition~1 we use induced wheels of even order and induced cycles of even order in a given embedded graph $G$ to describe small CDCs in the dual of $G$ or to show that there are many not-too-large CDCs. In~\cite{HS24}, Hu\v{s}ek and \v{S}\'amal use induced wheels of $G$---that is, in the dual: a face and all faces adjacent to it---to show that there are many CDCs in the dual of $G$. Proposition~1 is more general, as the induced cycles we treat may be such that there is a vertex $v$ adjacent to each of the cycle's vertices, so that the vertex together with the cycle induce a wheel (in this case we would be exactly in the framework of~\cite{HS24}), but this is certainly not a requirement. We point out that the above version of our proposition is in fact conservative in the sense that it could be generalised by letting some intersections between cycles in the dual occur. However, this would make the statement substantially more technical and cumbersome without much benefit to the following two applications we wish to here focus on.


\bigskip

\noindent \textbf{Theorem 4.} \emph{Every planar $2$-connected cubic graph on $n > 4$ vertices has a CDC of size at most~$\frac{n}{2}$.}

\begin{proof} We first need an auxiliary result on joining CDCs.

\smallskip

\noindent \textsc{Claim.} \textit{Let $T_1$ and $T_2$ be disjoint plane triangulations and $\Delta_1$ ($\Delta_2$) a facial triangle in $T_1$ ($T_2$). We can identify $\Delta_1$ and $\Delta_2$ such that a plane triangulation $T$ is obtained. For $i \in \{ 1, 2 \}$, if the cubic graph $T_i^*$ admits a $t_i$-CDC ${\cal C}_i$, then $T^*$ admits a $(t_1 + t_2 - 3)$-CDC ${\cal C}'$.}

\smallskip

\noindent \textit{Proof of the Claim.} Throughout the proof, $i \in \{ 1, 2 \}$. Let $x_i$ be the vertex in $T_i^*$ corresponding to $\Delta_i$. There are precisely three pairwise distinct cycles $C_i^1, C_i^2, C_i^3$ in ${\cal C}_i$ which cover twice the three edges incident to $x_i$. We now consider $T_i$ to be a subgraph of $T$. For every $j \in \{ 1, 2, 3 \}$ we remove from $C^j_i$ the vertex $x_i$ and join $C^j_1 - x_1$ with $C^k_2 - x_2$ in the obvious way, for an appropriate $k \in \{ 1, 2, 3 \}$ (this is possible as there is essentially only one way in which the three cycles can be distributed among the three edges incident with $x_i$). This gives three new cycles $C^j$. The cycles
$$\left( \bigcup_{i \in \{ 1, 2 \}} {\cal C}_i \setminus \{ C_i^1, C_i^2, C_i^3 \} \right) \cup \{ C^1, C^2, C^3 \}$$
are a CDC of $T^*$. Since we joined, in a bijective manner, three cycles in ${\cal C}_1$ with three cycles in ${\cal C}_2$, its size is $t_1 + t_2 - 3$.


\smallskip

In the remainder of this proof, all graphs are considered to be planar graphs embedded in the plane. Let $G$ be a 2-connected $n$-vertex cubic simple graph which is not $K_4$ and $G^*$ its dual, which is a (possibly non-simple) triangulation. The graph $G^*$ is both simple and 3-connected iff $G$ is 3-connected. 


We first assume $G^*$ to be 4-connected. Let $v$ be a vertex in $G^*$. If the degree of $v$ is even, we consider $G^*[N(v)]$, a cycle of even length which must be chordless as $G^*$ is 4-connected. If the degree of $v$ is odd, we consider the subgraph $W$ of $G^*$ given by $G^*[N[v]]$. The graph $W$ has even order at least 6, and as $G^*$ is 4-connected, it is an induced wheel. So regardless of whether the degree of $v$ is even or not, by Proposition~1 the graph $G$ contains a $\frac{n}{2}^-$-CDC. We will summarise this paragraph's findings by ($\dagger$).

Below we will use the same decomposition argument for plane triangulations as in Section~2.2. Let $T$ be a plane triangulation. Then either $T$ is 4-connected or there exists a separating triangle $S$ in $T$ such that $S$ together with its interior is either $K_4$ or a plane 4-connected triangulation.

We have proven that (${\frak B}$) holds if $G^*$ is 4-connected, so it remains to deal with the cases when $G^*$ has connectivity lower than 4. Assume $G^*$ has connectivity~3. If $G^*$ contains at most one cubic vertex, consider a separating triangle $S$ in $G^*$ such that, without loss of generality, $S$ together with its interior, i.e.~${\rm Int}(S)$, is 4-connected (we have ${\rm Int}(S) \ne K_4$ since we require $G^*$ to contain at most one cubic vertex). Put $s := |F({\rm Int}(S))|$. By ($\dagger$), the planar cubic graph $({\rm Int}(S))^*$ admits a $\frac{s}{2}^-$-CDC. The plane cubic graph $({\rm Ext}(S))^*$ has, by Euler's formula, a $(\frac{s'}{2} + 2)^-$-CDC (simply by considering the cycles constituting face boundaries), where $s' := |F({\rm Ext}(S))|$. By the above Claim, $G$ then has a CDC of size at most $\frac{s}{2} + \frac{s'}{2} - 1 = \frac{n}{2}$, as desired. 

So we may suppose that $G^*$ contains at least two cubic vertices. We first observe that these must be non-adjacent as $G$ is a cubic graph on at least six vertices and thus at least five faces, so $G^*$ is a triangulation on at least five vertices. Hence, $G$ contains two disjoint triangles. By Euler's formula the face boundaries in $G$ yield an $(\frac{n}{2} + 2)$-CDC ${\cal C}$ of $G$. Removing from ${\cal C}$ the two cycles corresponding to the boundary cycles of these two triangles and rerouting the boundary cycles of the three adjacent faces in a straightforward way (each visits the two edges of the triangle it did not visit before), one can easily modify ${\cal C}$ into an $\frac{n}{2}$-CDC of $G$; as $G$ has no bridge, the three faces surrounding a triangle must be pairwise distinct.


Henceforth, we assume that $G$ has a 2-vertex-cut. We observe that Bondy had shown in~\cite{Bo90} that if (${\frak B}$) holds for all simple 3-connected 3-regular graphs, then it holds for all simple 2-connected 3-regular graphs. For the connectivity~2 case we will use an argument similar to the one given by Bondy, but include it nonetheless. That $G$ has a 2-vertex-cut implies, as $G$ is cubic, that $G$ has a 2-edge-cut $M$. Above we have already argued that if $G$ contains two disjoint triangles, we obtain the desired cover, so we may assume that there exists a component of $G - M$ such that none of its bounded faces is a triangle. Moreover, by elementary reduction arguments we may assume that this component, which we call $K$, is 2-connected, and that if $v$ and $w$ are the (necessarily distinct) vertices in $K$ incident with edges from $M$, we may suppose $v$ and $w$ to be non-adjacent; we see $K$ as a plane 2-connected graph of maximum degree~3. 







We consider the cubic plane 3-connected graph $H$ obtained by taking $K$ and joining its only 2-valent vertices, namely $v$ and $w$, by an edge. The graph $H$ has at least two bounded faces, at least one of which is not a triangle, so $H \ne K_4$. By the above result stating that every planar 3-connected cubic graph on $n$ vertices which is not $K_4$ has a CDC of size at most $\frac{n}{2}$, $H$ has a CDC ${\cal C}'$ of size at most $|V(H)|/2$. 

We now describe how to extend ${\cal C}'$ in order to obtain an $\frac{n}{2}^-$-CDC of $G$. Let $K' \ne K$ be the other component of $G - M$. As $K'$ is a subgraph of $G$, the embedding of $G$ into the plane gives an embedding of $K'$ in the plane; henceforth we see $K'$ as being this specific embedding of the corresponding abstract graph. Let $B$ be the boundary cycle of the infinite face of $K'$. The end-vertices $v', w'$ of $M$ which are not $v$ and $w$ edge-partition $B$ into $v'w'$-paths $B^1$ and $B^2$, with $v'$ ($w'$) adjacent to $v$ ($w$). In $H$, there are exactly two distinct cycles $C^1, C^2 \in {\cal C}'$ containing the edge $vw$. In $G$, for $i \in \{ 1, 2 \}$ we add to the $vw$-path $C^i - vw$ the edges $ww'$ and $vv'$ as well as the $v'w'$-path $B^i$, thus obtaining two cycles $D^1, D^2$. We add to ${\cal C}' \setminus \{ C^1, C^2 \} \cup \{ D^1, D^2 \}$ the boundary cycles of bounded faces of $K'$ to obtain the collection ${\cal C}''$ of cycles in $G$. Clearly, ${\cal C}''$ is a CDC. By Euler's formula its size is at most
$$|{\cal C}'| - 2 + 2 + |F(K')| - 1 = \frac{|V(H)|}{2} + \frac{n - |V(H)| + 2}{2} - 1 = \frac{n}{2},$$
which completes the proof. \end{proof}


The bound in Theorem~4 cannot be improved as noted by Bondy~\cite{Bo90}, see Fig.~4 for his example. 

\begin{center}
\includegraphics[height=25mm]{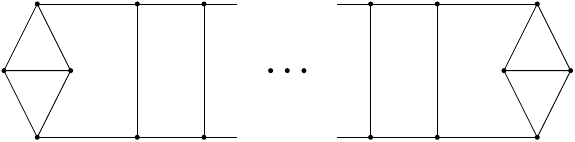}\\
Fig.~4: Infinitely many planar 2-connected cubic graphs in which every CDC has size at least~$\frac{n}{2}$.\\
As mentioned at the beginning of this section, a CDC in a cubic graph contains at most $\frac{n}{2} + 2$ cycles, so every CDC in the depicted family of graphs has size either~$\frac{n}{2}$, or~$\frac{n}{2}+1$, or~$\frac{n}{2}+2$.
\end{center}







\bigskip

Corollary 20 from~\cite{HS24} by Hu\v{s}ek and \v{S}\'amal states that given a fixed surface $\Sigma$, every bridgeless cubic graph embedded on $\Sigma$ and with representativity at least 4 has $2^{\Omega(\sqrt{n})}$ CDCs. We now use Proposition~1 to strengthen this result in the following way: we prove the same asymptotic behaviour under the additional assumption that the CDCs of the graph $G$ in question are of size at most the number of faces of $G$. We only sketch the proof as the arguments are very similar to the ones given by Hu\v{s}ek and \v{S}\'amal in~\cite{HS24}, combined with our Proposition~1.

In fact, we will only require a restricted version of Proposition~1, in the following sense. Given an induced wheel $W$ in an embedded graph $G$, we denote its central vertex by $v$. Then the face in $G^*$, the dual of $G$, corresponding to $v$ will be called \textit{central}. Similarly, given an induced cycle $C$ in an embedded graph $G$ such that every vertex in $C$ is adjacent to some vertex $v$, the face in $G^*$ corresponding to $v$ will be called \textit{central}. Note that the induced cycles in Proposition~1 do not have this adjacency restriction, but subsequent arguments will be easier to formulate if we impose it; it is also clear that these special cycles we are now treating are obtained by considering an induced wheel and removing its central vertex.

Like Hu\v{s}ek and \v{S}\'amal, we consider the vertex colouring of the square of the dual graph; we are interested in the graph's square in order to guarantee the disjointness of the aforementioned induced wheels and induced cycles. Amini, Esperet, and van den Heuvel~\cite{AEV13} proved that given a fixed surface $\Sigma$, there exists a positive real $c$ such that any graph $G$ embedded into $\Sigma$ satisfies $\chi(G^2) \le c\Delta(G)$. So there is a colour class in $G^2$ of size at least $|V(G)|/(c\Delta(G))$. The faces in $G^*$ corresponding to the vertices in this colour class can be chosen as being central for suitable wheels and cycles in $G$. Put $n := |V(G)|$. Again we follow Hu\v{s}ek and \v{S}\'amal and treat the following two cases: 

If $\Delta \ge \sqrt{n}$, we choose $F$ of size $\Delta$ to be a central face of a suitable induced wheel or cycle, depending on the parity of $\Delta$. Hu\v{s}ek and \v{S}\'amal established (see~\cite{HS24}, \cite{HS}) that the graph then contains $2^{\Omega(\Delta)}$ CDCs, i.e.\ in this case $2^{\Omega(\sqrt{n})}$ CDCs. It follows from the method given by Hu\v{s}ek and \v{S}\'amal that these CDCs are of size not greater than the number of faces of $G$, i.e.\  $\frac{n}{2} - 2\gamma + 2$.

If $\Delta \le \sqrt{n}$ we can find a colour class of size $\lceil |V(G)|/(c\Delta(G)) \rceil \ge \lceil n/(c\sqrt{n}) \rceil = \lceil \sqrt{n}/c \rceil$. We see the faces corresponding to the vertices in this colour class as central faces. We now consider the induced wheels and induced cycles around these faces, of which there are at least $\sqrt{n}/c$. Using the notation from Proposition~1, we have that $k + m \ge \sqrt{n}/c$ for some positive constant $c$. Thus, we have shown the following result.

\bigskip

\noindent \textbf{Proposition 2.} \emph{If the Strong Embedding Conjecture holds, then every $n$-vertex bridgeless cubic graph of genus $\gamma$ has $2^{\Omega(\sqrt{n})}$ CDCs of size at most $\frac{n}{2} - 2\gamma + 2$.}

\bigskip

We conclude the article with an observation on another question of Bondy.

\subsection{On the size of CDCs of \boldmath{3}-connected cubic graphs}

A natural follow-up question to conjecture (${\frak B}$) posed by Bondy himself is Remark 5.4 in his paper~\cite{Bo90}. Therein, he writes that he knows of no cubic 3-connected graph on $n$ vertices which fails to admit a CDC ${\cal C}$ such that $|{\cal C}| \le (n + 10)/4$. We were not able to prove or disprove this conjecture---what we could do is show that one cannot replace $(n + 10)/4$ by a smaller quantity. We first prove a useful lemma. Let $G$ and $H$ be cubic graphs, $x \in V(G)$ and $y \in V(H)$. We denote by $(G - x) \equiv (H - y)$ the graph resulting from connecting the three dangling edges of $G - x$ with the three dangling edges of $H - y$. When we do not care about which vertices are chosen, we simply write $G \equiv H$ for an arbitrarily chosen graph from all possible graphs obtained in this way.

\bigskip

\noindent \textbf{Lemma 3.} \textit{Let $G$ and $H$ be disjoint cubic graphs. Then $c(G \equiv H) = c(G) + c(H) - 3$.}

\begin{proof} 
Put $\Gamma := G \equiv H$. Let $x$ ($y$) be the vertex in $G$ ($H$) that we removed to construct $\Gamma$, and $M$ the 3-edge-cut between $G - x$ and $H - y$. To see that $c(\Gamma) \le c(G) + c(H) - 3$, we take a minimum cardinality CDC of $G$ and a minimum cardinality CDC of $H$, and merge them in the obvious way to get a CDC of $\Gamma$. For $c(\Gamma) \ge c(G) + c(H) - 3$, note that in any CDC of $\Gamma$ there are three pairwise disjoint sets of cycles: a) cycles contained entirely within one half, b) cycles contained entirely within the other half and c) precisely three cycles containing two of the edges of the 3-edge-cut $M$; when we say ``half'', we refer to the graph $\Gamma$ where all old vertices that belonged to $G$ (or $H$) are removed. Now indeed, the size of the set containing the cycles of one half cannot be smaller than $c(G)-3$ (or $c(H)-3$), because then we could construct a CDC of $G$ (or $H$) containing fewer than $c(G)$ (or $c(H)$) cycles by contracting the other half to a single vertex.
\end{proof}

Applying iteratively Lemma~3 to copies of the Petersen graph, in which the smallest CDC has size~5, we obtain:

\bigskip

\noindent \textbf{Proposition 3.} \emph{There exist infinitely many cubic $3$-connected graphs $G$ in which every CDC has size at least $(|V(G)|+10)/4.$}

\section{Notes}

\noindent \textbf{1.} Regarding the construction we used for Proposition~3, Bondy's remarks suggest that he knew of such a family; we include it here for completeness. 
He writes in~\cite{Bo90} that a ``still stronger conjecture may well be valid for 3-connected cyclically 4-edge-connected 3-regular graphs.'' Our small-scale computations did not yield any counterexamples to the $(n+10)/4$ bound, and likewise did not provide us with an inspiration for what the right bound for the 3-connected cyclically 4-edge-connected 3-regular case could be.

\smallskip

\noindent \textbf{2.} Hu\v{s}ek and \v{S}\'amal~\cite{HS24} conjecture that every bridgeless cubic graph with $n$ vertices has at least $2^{(n/2)-1}$ CDCs. We confirmed the computations of Hu\v{s}ek and \v{S}\'amal yielding that there are no small counterexamples, and also attacked the problem by relaxing the regularity condition. It turns out that it is easy to find counterexamples if one only imposes that the graphs have minimum degree~3. However, they all have connectivity~1. So we formulate here, as a question, the following variation of the aforementioned conjecture: Does every 2-connected graph on $n$ vertices and with minimum degree at least 3 have at least $2^{(n/2)-1}$ CDCs?

\smallskip

\noindent \textbf{3.} Next to questions of Bondy, Seyffarth's theorem was the main motivation for this paper. It remains open whether the ``$n-1$'' in its statement can be improved to ``$n-2$'' (and this remains open even when restricting the problem to triangulations), which would be best possible. But what is the right bound for the \textit{true} version of Seyffarth's theorem, i.e.\ when dealing only with CDCs whose cycles are pairwise distinct? Clearly, every plane 4-connected graph with $f$ faces contains a true $f$-CDC. For antiprisms of order $n$, we have $f = n + 2$, and there are no planar graphs of minimum degree at least 4 with $f < n + 2$ by Euler's formula. We wonder whether perhaps every planar 4-connected graph contains a \textit{true} $(n+2)^-$-CDC. 

\smallskip

\noindent \textbf{4.} Related to the previous remark, although Corollary~2 does not provide \textit{substantial} evidence, we do feel compelled to raise the intriguing question which other plane graphs with $f$ faces might have a unique true $f$-CDC (which then must be exactly the set of face boundaries). On the one hand, it is obviously not true for all planar 2-connected graphs, and it is not difficult to describe counterexamples for the connectivity~3 case. Likewise, toroidal 4-connected counterexamples are easily found. On the other hand, every 4-connected triangulation with $f$ faces has a unique true $f$-CDC, as a second true $f$-CDC would have to contain a cycle of length greater than~3 since there are no separating triangles. Furthermore, all graphs described in Theorem~2~(i)---which are close to triangulations in the sense that all but two faces are triangular---exhibit the same behaviour. But we could not determine whether every plane 4-connected graph with $f$ faces has a unique true $f$-CDC. 

\smallskip

\noindent \textbf{5.} Finally, we deliver on our promise made at the beginning of Section~2 to give a short proof of the fact that the edges of every Hamiltonian graph on at least three vertices can be covered by three even subgraphs. Here we allow loops and multiple edges to occur. Suppose $G$ is a counterexample having the least size. Let ${\frak h}$ be a Hamiltonian cycle of $G$. Clearly $G$ cannot have a vertex of degree~1, and $G$ cannot have a vertex $v$ of degree~2 with incident edges $uv$ and $vw$ because then $(V(G) \setminus \{ v \}, E(G) \setminus \{ uv,vw \} \cup \{uw\})$ would be a smaller counterexample. Similarly, $G$ cannot have a vertex of degree at least~4 because then we could remove two edges incident with this vertex which are not contained in ${\frak h}$ and replace them by a single edge connecting their two other endpoints, and we would again have a smaller counterexample. Thus $G$ is cubic and we can colour the edges of a Hamiltonian cycle with colours 1 and 2, and all other edges with colour 3. Taking the three pairs of colours, we obtain three even subgraphs covering every edge of $G$ exactly twice, a contradiction.

\bigskip

\bigskip

\noindent\textbf{Acknowledgements.} Jorik Jooken is supported by an FWO grant with grant number 1222524N. 
 Ben Seamone is supported by grants from NSERC (RGPIN-2017-06673) and FRQNT (CO-289707).  The authors thank On-Hei Solomon Lo, who participated in fruitful preliminary discussions on the manuscript's topics.

\end{document}